\documentclass{ifacconf}

\usepackage{enumitem}
\usepackage{siunitx}
\usepackage{tikz}
\usepackage[siunitx]{circuitikz}
\usepackage{comment}
\usepackage[tbtags]{amsmath}
\usepackage[tbtags]{mathtools}
\usepackage{tabularx} 
\usepackage{diagbox}
\usepackage{graphicx}
\usepackage{mathrsfs}
\usepackage{cite}
\usepackage{amssymb,amsfonts}
\mathtoolsset{showonlyrefs=true}

\usepackage{mathdots}
\usepackage{algorithm}
\usepackage{algpseudocode}
\algrenewcommand\algorithmicrequire{\textbf{Given:}}
\algrenewcommand\algorithmicensure{\textbf{Return:}}
\usepackage{float}
\usepackage{accents}
\usepackage{multicol}
\usepackage{multirow}
\usepackage{cite}

\usepackage{array,multirow}
\usepackage{hhline}
\usepackage{arydshln}

\newtheorem{theorem}{Theorem}
\newtheorem{corollary}[theorem]{Corollary}
\newtheorem{lemma}[theorem]{Lemma}

\newtheorem{proposition}[theorem]{Proposition}
\newtheorem{problem}{Problem}

\newtheorem{definition}[theorem]{Definition}

\DeclareMathOperator{\rank}{rank}

\DeclareMathOperator{\im}{im}

\newcommand{\set}[2]{\left\{#1 \mid #2\right\}}

\usepackage{graphicx}      
\usepackage{natbib}        
\begin{document}
\begin{frontmatter}

\title{Experiment design using prior knowledge on controllability and stabilizability} 

\thanks[second]{The work of Henk van Waarde was supported by
the Dutch Research Council under the NWO Talent Programme Veni
Agreement (VI.Veni.22.335).}

\author[First]{Amir Shakouri} 
\author[First,second]{Henk J. van Waarde} 
\author[First]{M. Kanat Camlibel}

\address[First]{Bernoulli Institute for Mathematics, Computer Science and Artificial Intelligence, University of Groningen, The Netherlands (e-mail: {\tt a.shakouri@rug.nl}, {\tt h.j.van.waarde@rug.nl}, {\tt m.k.camlibel@rug.nl}).}

\begin{abstract}
In this paper, we consider the problem of designing input signals for an unknown linear time-invariant system in such a way that the resulting input-state data is suitable for identification or stabilization. We will take into account prior knowledge on system-theoretic properties of the system, in particular, controllability and stabilizability. For this, we extend the notion of \emph{universal inputs} to incorporate prior knowledge on the system. An input is called universal for identification (resp., stabilization) if, when applied to any system complying with the prior knowledge, it results in data suitable for identification (resp., stabilization) regardless of the initial condition. We provide a full characterization of such universal inputs. In addition, we discuss online experiment design using prior knowledge, and we study cases where this approach results in the shortest possible experiment for identification and stabilization. 
\end{abstract}

\begin{keyword}
Experiment design, system identification, data-driven control, prior knowledge.
\end{keyword}

\end{frontmatter}

\section{Introduction}

Data collected from an unknown linear time-invariant system can be used to identify the dynamics, or to directly analyze system properties and design controllers, see \citep{markovsky2021behavioral,van2025data,berberich2025overview} and the references therein. The feasibility of such data-driven methods depends on the richness of the data. If the data does not contain sufficient information about the system, data-driven modeling and control are not possible. This motivates studying the \emph{experiment design} problem that answers the following question: How to select input signals in such a way that the resulting data is suitable for modeling and control?

For controllable systems, \cite{willems2005note} provides an offline method for the experiment design for system identification. Due to Willems et al.'s \emph{fundamental lemma}, a finite number of input-output samples, generated by a persistently exciting input signal, can be used to parametrize all finite-length trajectories of a controllable system. As a corollary, a persistently exciting input of a certain order guarantees that the dynamics of a controllable system can be uniquely identified from the resulting input-output data. An interesting property of persistently exciting inputs is that they are \emph{universal} in the sense that they simultaneously work for \emph{all} controllable systems. Up until recently, a full characterization of universal inputs had been missing from the literature. This gap was filled in \citep{shakouri2025new}, where it is shown that an input is universal \emph{if and only if} it is persistently exciting of a sufficiently high order. 

An offline experiment design has the advantage of being simple and computationally inexpensive due to its open-loop nature. However, such advantages come with the cost of requiring a rather large number of data samples. In contrast, an \emph{online} experiment design selects the input at each step based on the collected data from the previous steps. Since this past data contains information about the dynamics, online approaches tailor the experiment design to the specific (unknown) system and to the specific state in which the system is initialized. As a consequence, online experiment design often outperforms offline methods in terms of sample efficiency. The first method in this line of work was developed by \cite{van2021beyond}. For controllable systems, the online experiment design algorithm presented in \citep{van2021beyond} generates the \emph{shortest} possible \emph{input-state} data suitable for system identification. It has also been shown that an online experiment design can generate the shortest experiment for system identification using \emph{input-output} data, for which the reader can refer to \citep{van2021beyond,camlibel2024beyond,camlibel2025shortest}. 

It is well-known that if the system is not controllable, then collecting input-output data suitable for system identification is not possible in general. In fact, this depends on the initial state of the system, see \cite[Thm. 1]{yu2021controllability}. Nevertheless, it was shown recently by \cite{gramlich2023fast} that if the system is known to be stabilizable, then a persistently exciting input guarantees that the collected input-state data can be used to find a stabilizing feedback for the unknown system. It was also shown in the same work that an online experiment design can accomplish the same goal. Therefore, \emph{prior knowledge} on the true system, such as controllability or stabilizability, plays an important role in the experiment design for identification or stabilization. So far, a characterization of universal inputs using such system-theoretic prior knowledge is missing from the literature. Moreover, it is not yet known how to design online experiments that yield the shortest dataset required for stabilization.

In this work, we extend the notion of universal inputs in \citep{shakouri2025new} to take the prior knowledge into account. An input is called universal for identification (resp., stabilization) if, when applied to any system complying with the prior knowledge, it results in data suitable for identification (resp. stabilization) regardless of the initial condition. In particular, we consider prior knowledge of controllability and stabilizability. We prove that an input is universal for identification (resp. stabilization) using controllability (resp. stabilizability) as prior knowledge \emph{if and only if} it is persistently exciting of a certain order. In addition, we study online experiment design methods using such prior knowledge. We show that if the initial state satisfies a certain condition, then an online experiment design method generates the \emph{shortest} dataset required for \emph{stabilization}. For this, we leverage the recently developed necessary and sufficient conditions for data-driven stabilization using prior knowledge on controllability and stabilizability in \citep{shakouri2025pk}.

\section{Preliminaries}
\label{sec:II}

\subsection{Notation and basic terminology}

Let $\mathbb{Z}_+$ and $\mathbb{N}$ denote the sets of nonnegative and positive integers, respectively. For $i,j\in\mathbb{N}$ with $i\leq j$, we define the integer interval $[i,j]\coloneqq\set{k\in\mathbb{Z}}{i\leq k\leq j}$. Given $v:[i,j]\rightarrow\mathbb{R}^n$, we define $V_{[i,j]}\coloneqq\begin{bmatrix}
v(i) & v(i+1) & \cdots & v(j)
\end{bmatrix}\in\mathbb{R}^{n\times (j-i+1)}$ and
$v_{[i,j]}\coloneqq\begin{bmatrix}
v(i)^\top & v(i+1)^\top & \cdots & v(j)^\top
\end{bmatrix}^\top\in\mathbb{R}^{n(j-i+1)}$.
For $T\in\mathbb{N}$, we say that $v_{[0,T-1]}$ is \emph{persistently exciting of order $k$} if the Hankel matrix of $v_{[0,T-1]}$ of depth $k$, $\begin{bmatrix}
v_{[0,k-1]} & v_{[1,k]} & \cdots & v_{[T-k,T-1]}
\end{bmatrix}$, has full row rank.

Let $n,m\in\mathbb{N}$. Consider the input-state system
\begin{equation}
\label{eq:1}
x(t+1)=Ax(t)+Bu(t),
\end{equation}
where $t\in\mathbb{Z}_+$ denotes time, $u(t)\in\mathbb{R}^m$ is the input, and $x(t)\in\mathbb{R}^n$ is the state. We identify system \eqref{eq:1} with the pair of matrices $(A,B)\in\mathcal{M}$, where
\begin{equation*}
\mathcal{M}\coloneqq \mathbb{R}^{n\times n}\times\mathbb{R}^{n\times m}.
\end{equation*}
We denote the reachable subspace of $(A,B)\in\mathcal{M}$ by $\mathcal{R}(A,B)\coloneqq\im\begin{bmatrix}
B & AB & \cdots & A^{n-1}B
\end{bmatrix}$. We also denote the sets of controllable and stabilizable systems in $\mathcal{M}$, respectively, by $\Sigma_\text{cont}$ and $\Sigma_\text{stab}$. 

We denote the \emph{input-state behavior} of \eqref{eq:1} by
\begin{equation}
\mathfrak{B}(A,\!B)\!\coloneqq\!
\set{(u,x)\!:\!\mathbb{Z}_+\!\rightarrow\! \mathbb{R}^m\!\times\!\mathbb{R}^n}{\!\eqref{eq:1}\text{ holds for all }t\!\in\!\mathbb{Z}_+\!}\!,
\end{equation}
and its \emph{$k$-restricted behavior} by
\begin{equation}
\mathfrak{B}_k(A,B)\coloneqq\set{(u_{[0,k-1]} , x_{[0,k]})}{(u,x)\in\mathfrak{B}(A,B)}.
\end{equation}

\subsection{Data informativity and universal inputs}

We consider the \emph{true system},
\begin{equation}
(A_\text{true},B_\text{true})\in\mathcal{M},
\end{equation}
which is assumed to be unknown. We collect input-state data of length $T\in\mathbb{N}$ from the true system, denoted by
\begin{equation}
\mathcal{D}=(u_{[0,T-1]} , x_{[0,T]})\in\mathfrak{B}_T(A_\text{true},B_\text{true}).
\end{equation}

Given $\mathcal{D}$, we define the set of \emph{data-consistent systems} as
\begin{equation}
\Sigma_\mathcal{D}\coloneqq\set{(A,B)\in\mathcal{M}}{\mathcal{D}\in\mathfrak{B}_T(A,B)}.
\end{equation}
Clearly, we have $(A_\text{true},B_\text{true})\in\Sigma_\mathcal{D}$. Moreover, we assume that we have prior knowledge on the true system, i.e.,
\begin{equation}
(A_\text{true},B_\text{true})\in\Sigma_\text{pk},
\end{equation}
where $\Sigma_\text{pk}\subseteq\mathcal{M}$ is a given set. In case no prior knowledge is given, we set $\Sigma_\text{pk}=\mathcal{M}$. 

The data $\mathcal{D}$ and the prior knowledge set $\Sigma_\text{pk}$ can be used to identify the true system or to find a stabilizing feedback gain for it. For both of these goals, we exploit the fact that the true system satisfies
\begin{equation}
(A_\text{true},B_\text{true})\in\Sigma_\mathcal{D}\cap\Sigma_\text{pk}.
\end{equation}
In case the set $\Sigma_\mathcal{D}\cap\Sigma_\text{pk}$ is a singleton, one can uniquely identify the true system. Moreover, in case there exists a single feedback gain that simultaneously stabilizes all systems within $\Sigma_\mathcal{D}\cap\Sigma_\text{pk}$, then such a feedback also stabilizes the true system. We formalize these concepts by extending the notions of data informativity in \citep{van2020data} to incorporate prior knowledge on the true system. 
\begin{definition}[Data informativity]
\label{def:inf}
Let $T\in\mathbb{N}$. Then, $\mathcal{D}$ is called
\begin{enumerate}[label=(\roman*),ref=\ref{def:inf}(\roman*)]
    \item\label{def:inf(i)} $\Sigma_\text{pk}$--\emph{informative for identification} if $\Sigma_\mathcal{D}\cap\Sigma_\text{pk}$ is a singleton. 
    \item\label{def:inf(ii)} $\Sigma_\text{pk}$--\emph{informative for stabilization} if there exists a \linebreak $K\in\mathbb{R}^{m\times n}$ such that $A+BK$ is Schur for all $(A,B)\in\Sigma_\mathcal{D}\cap\Sigma_\text{pk}$. 
\end{enumerate}
\end{definition}

We note that $\mathcal{M}$--informativity coincides with the informativity notions studied in \citep{van2020data,van2025data}.

The experiment design problem concerns finding an input signal by which the true system generates informative data. Since the true system is unknown, one has to design such an input so that it works for all systems within the prior knowledge set. To this end, we investigate inputs with the property that, when applied to any system in $\Sigma_\text{pk}$, the generated data is $\Sigma_\text{pk}$--informative for identification or stabilization regardless of the initial condition. To formalize this, we introduce the notion of universal inputs as follows. 

\begin{definition}[Universal inputs]
\label{def:uni}
Let $T\in\mathbb{N}$. An input $u_{[0,T-1]}$ is called $\Sigma_\text{pk}$--\emph{universal for identification} (resp., \emph{stabilization}) if for every $(A,B)\in\Sigma_\text{pk}$ and every $x_{[0,T]}$ satisfying $\mathcal{D}=(u_{[0,T-1]} , x_{[0,T]})\in\mathfrak{B}_T(A,B)$ we have that $\mathcal{D}$ is $\Sigma_\text{pk}$--informative for identification (resp., stabilization).
\end{definition}

\subsection{Problem statement}

An interesting problem that has not been addressed in the literature is to characterize $\Sigma_\text{pk}$--universal inputs for a given set of prior knowledge $\Sigma_\text{pk}$.

\begin{problem}
\label{prob:1}
Find necessary and sufficient conditions under which an input is 
(i) $\Sigma_\text{pk}$--universal for identification; (ii) $\Sigma_\text{pk}$--universal for stabilization.
\end{problem}

In this paper, we are mainly interested in prior knowledge that represents the controllability or stabilizability of the true system. We thus investigate the solutions to Problem~\ref{prob:1} for two cases:
\begin{equation}
\Sigma_\text{pk}=\Sigma_\text{cont}\ \text{ and }\ \Sigma_\text{pk}=\Sigma_\text{stab}.
\end{equation}

\section{Universal Experiment Design for Identification}
\label{sec:univ_id}

In this section, we provide the solution to Problem~\ref{prob:1}(i) for $\Sigma_\text{pk}=\Sigma_\text{cont}$ and $\Sigma_\text{pk}=\Sigma_\text{stab}$. To that end, first, we present necessary and sufficient conditions for the informativity of the data for identification using prior knowledge. For the sake of simplicity, we use the following notation:
\begin{equation}
U_-\coloneqq U_{[0,T-1]},\ 
X_-\coloneqq X_{[0,T-1]},\ \text{ and }\ X_+\coloneqq X_{[1,T]}.
\end{equation}

\begin{theorem}
\label{th:inf_id_pk}
Suppose that $\Sigma_\text{pk}$ is open. Then, the following statements are \emph{equivalent}:
\begin{enumerate}[label=(\alph*),ref=\ref{th:inf_id_pk}(\alph*)]
    \item\label{th:inf_id_pk(a)} $\mathcal{D}$ is $\Sigma_\text{pk}$--informative for identification.
    \item\label{th:inf_id_pk(b)} $\mathcal{D}$ is $\mathcal{M}$--informative for identification.
    \item\label{th:inf_id_pk(c)} $\rank \begin{bmatrix}
    X_-^\top & U_-^\top
    \end{bmatrix}=n+m$.
\end{enumerate}
\end{theorem}
\begin{pf}
It follows from \cite[Prop. 6]{van2020data} that 
(b) and (c) are equivalent. (b)$\Rightarrow$(a): Suppose that (b) holds, i.e., $\Sigma_\mathcal{D}$ is a singleton. Since the true system belongs to both $\Sigma_\mathcal{D}$ and $\Sigma_\text{pk}$, we have that $\Sigma_\mathcal{D}\cap\Sigma_\text{pk}$ is a singleton. Therefore, the data are $\Sigma_\text{pk}$--informative for identification. (a)$\Rightarrow$(c): Assume that (a) holds. Let $A_0\in\mathbb{R}^{n\times n}$ and $B_0\in\mathbb{R}^{n\times m}$ be such that $A_0X_-+B_0U_-=0$. This implies that $(A_\text{true}+\alpha A_0,B_\text{true}+\alpha B_0)\in\Sigma_\mathcal{D}$ for all $\alpha\in\mathbb{R}$. Since $(A_\text{true},B_\text{true})\in\Sigma_\mathcal{D}\cap\Sigma_\text{pk}$ and $\Sigma_\text{pk}$ is an open set, there exists a sufficiently small $\varepsilon>0$ such that $(A_\text{true}+\alpha A_0,B_\text{true}+\alpha B_0)\in\Sigma_\mathcal{D}\cap\Sigma_\text{pk}$ for all $|\alpha|\leq\varepsilon$. Since $\Sigma_\mathcal{D}\cap\Sigma_\text{pk}$ is a singleton due to (a), we see that both $A_0$ and $B_0$ must be zero. Therefore, (c) holds. \hfill $\blacksquare$
\end{pf}

Since $\Sigma_\text{cont}$ and $\Sigma_\text{stab}$ are both open sets, the following corollary is an immediate consequence of Theorem~\ref{th:inf_id_pk}.

\begin{corollary}
\label{cor:inf_id}
The following statements hold:
\begin{enumerate}[label=(\alph*),ref=\ref{cor:inf_id}(\alph*)]
    \item\label{cor:inf_id(a)} Suppose that $(A_\text{true},B_\text{true})\in\Sigma_\text{stab}$. Then, $\mathcal{D}$ is $\Sigma_\text{stab}$--informative for identification \emph{if and only if} it is $\mathcal{M}$--informative for identification.
    \item\label{cor:inf_id(b)} Suppose that $(A_\text{true},B_\text{true})\in\Sigma_\text{cont}$. Then, $\mathcal{D}$ is $\Sigma_\text{cont}$--informative for identification \emph{if and only if} it is $\mathcal{M}$--informative for identification.
\end{enumerate}
\end{corollary}

Now, we turn our attention to universal experiment design for identification. First, we consider the case $\Sigma_\text{pk}=\Sigma_\text{cont}$. Based on Corollary~\ref{cor:inf_id}, a $\Sigma_\text{cont}$--universal input, $u_{[0,T-1]}$, has to guarantee that for every $(A,B)\in\Sigma_\text{cont}$ and every $x_{[0,T-1]}$ satisfying $\mathcal{D}=(u_{[0,T-1]},x_{[0,T]})\in\mathfrak{B}_T(A,B)$, we have that $\mathcal{D}$ is $\mathcal{M}$--informative for identification. A sufficient condition for this can be inferred from a corollary of Willems et al.'s fundamental lemma, stated next.

\begin{lemma}[\cite{willems2005note}]
\label{lem:Willems}
Let $(A,B)\in\Sigma_\text{cont}$ and $u_{[0,T-1]}$ be persistently exciting of order $n+1$. Then, for every $x_{[0,T]}$ satisfying $(u_{[0,T-1]},x_{[0,T]})\in\mathfrak{B}_T(A,B)$ we have $\rank \begin{bmatrix}X_-^\top & U_-^\top\end{bmatrix}=n+m$.
\end{lemma}

It follows from Corollary~\ref{cor:inf_id} and Lemma~\ref{lem:Willems} that persistently exciting inputs of order $n+1$ are $\Sigma_\text{cont}$--universal for identification. The fundamental lemma, however, does not provide a full characterization of such universal inputs. Such a characterization was provided recently by \cite{shakouri2025new}, where a converse to Willems' fundamental lemma was presented as follows. 

\begin{lemma}[\cite{shakouri2025new}]
\label{lem:shakouri}
If $u_{[0,T-1]}$ is \emph{not} persistently exciting of order $n+1$. then there exist $(A,B)\in\Sigma_\text{cont}$ and $x_{[0,T]}$ such that $(u_{[0,T-1]},x_{[0,T]})\in\mathfrak{B}_T(A,B)$ and $\rank X_-<n$.
\end{lemma}

Lemmas~\ref{lem:Willems} and~\ref{lem:shakouri} now lead to a full characterization of $\Sigma_\text{cont}$--universal inputs presented next. 

\begin{theorem}
\label{th:1}
Let $T\in\mathbb{N}$ and $u_{[0,T-1]}\in\mathbb{R}^{Tm}$. Then, the following statements are \emph{equivalent}:
\begin{enumerate}[label=(\alph*),ref=\ref{th:1}(\alph*)]
    \item\label{th:1(a)} $u_{[0,T-1]}$ is \mbox{$\Sigma_\text{cont}$--universal} for identification.
    \item\label{th:1(b)} $u_{[0,T-1]}$ is persistently exciting of order $n+1$. 
\end{enumerate}
\end{theorem}
\begin{pf}
We first prove that (b) implies (a). Suppose that (b) holds. Let $(A,B)\in\Sigma_\text{cont}$ and $x_{[0,T-1]}$ be such that $\mathcal{D}=(u_{[0,T-1]},x_{[0,T]})\in\mathfrak{B}_T(A,B)$. It follows from Lemma~\ref{lem:Willems} that $\rank \begin{bmatrix} X_-^\top & U_-^\top \end{bmatrix}=n+m$. Due to Theorem~\ref{th:inf_id_pk}, $\mathcal{D}$ is $\mathcal{M}$--informative for identification. Hence, it follows from Corollary~\ref{cor:inf_id} that $\mathcal{D}$ is $\Sigma_\text{cont}$--informative for identification. Since this argument holds for all $(A,B)\in\Sigma_\text{cont}$ and all $x_{[0,T-1]}$ satisfying $(u_{[0,T-1]},x_{[0,T]})\in\mathfrak{B}_T(A,B)$, the input $u_{[0,T-1]}$ is $\Sigma_\text{cont}$--universal for identification. Now, we show that (a) implies (b), for which we use a proof by contraposition. Suppose that (b) does not hold. It follows from Lemma~\ref{lem:shakouri} that there exist $(A,B)\in\Sigma_\text{cont}$ and $x_{[0,T]}$ satisfying $\mathcal{D}=(u_{[0,T-1]},x_{[0,T]})\in\mathfrak{B}_T(A,B)$ such that $X_-^\top$, and thus $\begin{bmatrix}X_-^\top & U_-^\top\end{bmatrix}$, does not have full column rank. Due to Theorem~\ref{th:inf_id_pk}, $\mathcal{D}$ is not $\mathcal{M}$--informative for identification. Hence, it follows from Corollary~\ref{cor:inf_id} that $\mathcal{D}$ is also not $\Sigma_\text{cont}$--informative for identification. Therefore, the input $u_{[0,T-1]}$ is not $\Sigma_\text{cont}$--universal for identification. \hfill $\blacksquare$
\end{pf}

Now, we consider the case $\Sigma_\text{pk}=\Sigma_\text{stab}$. To study this, we first present the following theorem.

\begin{theorem}
\label{th:uni_uncontrollable}
Suppose that $\Sigma_\text{pk}$ is open and satisfies $\Sigma_\text{pk}\not\subseteq\Sigma_\text{cont}$. Then, there are no $\Sigma_\text{pk}$--universal inputs for identification. 
\end{theorem}
\begin{pf}
Let $T\in\mathbb{N}$ and $u_{[0,T-1]}\in\mathbb{R}^{mT}$. We show that $u_{[0,T-1]}$ is not $\Sigma_\text{pk}$--universal for identification. Let $(A,B)\in\Sigma_\text{pk}$ be not controllable. Take $x(0)=x_0$ such that $\mathcal{R}(A,\begin{bmatrix}
B & x_0 
\end{bmatrix})\neq \mathbb{R}^n$. Since $\im X_-\subseteq \mathcal{R}(A,\begin{bmatrix}
B & x_0 
\end{bmatrix})$, this implies that $\rank X_-<n$, and hence, $\rank\begin{bmatrix}X_-^\top & U_-^\top\end{bmatrix}<n+m$. Thus, it follows from Theorem~\ref{th:inf_id_pk} that $\mathcal{D}$ is not $\Sigma_\text{pk}$--informative for identification. Therefore, $u_{[0,T-1]}$ is not $\Sigma_\text{pk}$--universal for identification.\hfill $\blacksquare$
\end{pf}

A noteworthy corollary of Theorem~\ref{th:uni_uncontrollable} is the following.

\begin{corollary}
\label{cor:cont_is_largest}
There exists a $\Sigma_\text{pk}$--universal input for identification \emph{if and only if} $\Sigma_\text{pk}\subseteq\Sigma_\text{cont}$.
\end{corollary}

Due to Corollary~\ref{cor:cont_is_largest}, $\Sigma_\text{cont}$ is the \emph{largest} set for which there exists a universal input for identification. Hence, the following corollary is an immediate consequence of this result.

\begin{corollary}
\label{cor:uni_stab}
There are no $\Sigma_\text{stab}$--universal inputs for identification. 
\end{corollary}

\section{Universal Experiment Design for Stabilization}
\label{sec:univ_stab}

In this section, we provide the solution to Problem~\ref{prob:1}(ii) for $\Sigma_\text{pk}=\Sigma_\text{cont}$ and $\Sigma_\text{pk}=\Sigma_\text{stab}$. To that end, first, we present necessary and sufficient conditions for the informativity of the data for stabilization. We recall that in case no prior knowledge exists, \mbox{$\Sigma_\text{pk}=\mathcal{M}$}, one can fully characterize the informativity of the data as follows. 

\begin{proposition}[\cite{van2020data}]
\label{prop:inf0}
The data $\mathcal{D}$ is $\mathcal{M}$--informative for stabilization \emph{if and only if} $X_-$ has full row rank and there exists a right inverse $X_-^\dagger$ of $X_-$ such that $X_+X_-^\dagger$ is Schur.
\end{proposition}

A necessary condition for $\mathcal{M}$--informativity for stabilization is $\rank X_- = n$. As a result, $\mathcal{M}$--informativity for stabilization requires the number of data samples to be larger than or equal to the state dimension, i.e., $T\geq n$. 

It was shown recently in \citep{shakouri2025pk}, that $\Sigma_\text{pk}$--informativity for stabilization can also be fully characterized for the cases $\Sigma_\text{pk}=\Sigma_\text{cont}$ and $\Sigma_\text{pk}=\Sigma_\text{stab}$. We summarize those results in the following proposition. 

\begin{proposition}[\cite{shakouri2025pk}]
\label{prop:inf_stab_pk}
The following \linebreak statements hold:
\begin{enumerate}[label=(\alph*),ref=\ref{prop:inf_stab_pk}(\alph*)]
    \item\label{prop:inf_stab_pk(a)} Suppose that $(A_\text{true},B_\text{true})\in\Sigma_\text{cont}$. Then, $\mathcal{D}$ is $\Sigma_\text{cont}$--informative for stabilization \emph{if and only if} it is $\mathcal{M}$--informative for stabilization. 
    \item\label{prop:inf_stab_pk(b)} Suppose that $(A_\text{true},B_\text{true})\in\Sigma_\text{stab}$ and $\rank X_-=n$. Then, $\mathcal{D}$ is $\Sigma_\text{stab}$--informative for stabilization \emph{if and only if} it is $\mathcal{M}$--informative for stabilization. 
    \item\label{prop:inf_stab_pk(c)} Suppose that $(A_\text{true},B_\text{true})\in\Sigma_\text{stab}$ and $\rank X_-<n$. Then, $\mathcal{D}$ is $\Sigma_\text{stab}$--informative for stabilization \emph{if and only if} the following conditions hold:
    \begin{equation}
    \label{eq:inf_stab_cond}
    \im X_+\subseteq\im X_-\ \text{ and }\  \im\begin{bmatrix}
    X_- \\ U_-
    \end{bmatrix}=\im X_- \times \mathbb{R}^m.
    \end{equation}
\end{enumerate}
\end{proposition}

It follows from statements (b) and (c) of Proposition~\ref{prop:inf_stab_pk} that the number of data samples for $\Sigma_\text{stab}$--informativity for stabilization must satisfy $T\geq\min\{n,m\}$. This is an enhanced lower bound compared to that of $\mathcal{M}$--informativity for stabilization, $T\geq n$, which is due to incorporating prior knowledge of stabilizability. 

Now, we turn our attention to generating informative data for stabilization using prior knowledge of systems' stabilizability or controllability. First, we state the following result, which lifts the controllability condition in the fundamental lemma, yielding a more general statement. 

\begin{lemma}
\label{lem:PE_stab}
Let $(A,B)\in\mathcal{M}$ and $(u_{[0,T-1]},x_{[0,T]})\in\mathfrak{B}_T(A,B)$ for some $T\in\mathbb{N}$. Consider the following statements:
\begin{enumerate}
    \item[(i)] $\im\begin{bmatrix}
    X_- \\ U_-
    \end{bmatrix}=\mathcal{R}(A,\begin{bmatrix}
    B & x(0)
    \end{bmatrix})\times \mathbb{R}^m$.
    \item[(ii)] $\im\begin{bmatrix}
    X_- \\ U_-
    \end{bmatrix}=\im X_-\times \mathbb{R}^m$ and $\im X_+\subseteq\im X_-$.
    \item[(iii)] $u_{[0,T-1]}$ is persistently exciting of order $n+1$.
\end{enumerate}
Then, statements (i) and (ii) are equivalent. Moreover,  if (iii) holds, then (i) and (ii) hold. 
\end{lemma}
\begin{pf}
Suppose that (i) holds. This obviously implies that $\im X_-= \mathcal{R}(A,\begin{bmatrix}
B & x(0)
\end{bmatrix})$. Thus, $\im X_-$ is $A$--invariant and contains $\im B$, i.e., $A\im X_-+\im B\subseteq \im X_-$. Since $X_+=AX_-+BU_-$, we have $\im X_+\subseteq A\im X_-+\im B\subseteq \im X_-$. Therefore, (i) implies (ii). Now, suppose that (ii) holds. The first condition in (ii) implies that $\im X_+=A\im X_-+\im B$. This, together with the second condition $\im X_+\subseteq \im X_-$, implies that $A\im X_-+\im B\subseteq \im X_-$, i.e., $\im X_-$ is $A$--invariant and contains $\im B$. Since we also have $x(0)\in\im X_-$, we thus have $\mathcal{R}(A,\begin{bmatrix}B & x_0\end{bmatrix})\subseteq \im X_-$. Note that $\im X_-\subseteq\mathcal{R}(A,\begin{bmatrix}B & x_0\end{bmatrix})$. Hence, it holds that $\im X_-= \mathcal{R}(A,\begin{bmatrix}B & x_0\end{bmatrix})$. Therefore, (ii) implies (i). For the rest, we refer the reader to \cite[Thm. 1]{yu2021controllability}, where it was proven that (iii) implies (i). \hfill $\blacksquare$
\end{pf}

Now, we use Proposition~\ref{prop:inf_stab_pk} along with Lemma~\ref{lem:PE_stab} to have a full characterization of universal inputs for stabilization. 

\begin{theorem}
\label{th:2}
Let $T\in\mathbb{N}$ and $u_{[0,T-1]}\in\mathbb{R}^{Tm}$. Then, the following statements are \emph{equivalent}:
\begin{enumerate}[label=(\alph*),ref=\ref{th:2}(\alph*)]
    \item\label{th:2(a)} $u_{[0,T-1]}$ is \mbox{$\Sigma_\text{stab}$--universal} for stabilization.
    \item\label{th:2(b)} $u_{[0,T-1]}$ is \mbox{$\Sigma_\text{cont}$--universal} for stabilization.
    \item\label{th:2(c)} $u_{[0,T-1]}$ is persistently exciting of order $n+1$. 
\end{enumerate}
\end{theorem}
\begin{pf}
(a)$\Rightarrow$(b): This follows immediately from the fact that $\Sigma_\text{cont}$ is a nonempty subset of $\Sigma_\text{stab}$. 

(b)$\Rightarrow$(c): To show this, we use a proof by contraposition. Assume that (c) does not hold, i.e., $u_{[0,T-1]}$ is not persistently exciting of order $n+1$. Due to Lemma~\ref{lem:shakouri}, there exist $(A,B)\in\Sigma_\text{cont}$ and $x_{[0,T-1]}$ satisfying $\mathcal{D}=(u_{[0,T-1]},x_{[0,T]})\in\mathfrak{B}_T(A,B)$ such that $\rank X_-<n$. Since $X_-$ does not have full row rank, it follows from Proposition~\ref{prop:inf0} that $\mathcal{D}$ is not $\mathcal{M}$--informative for stabilization. Now, it follows from Proposition~\ref{prop:inf_stab_pk(a)} that $\mathcal{D}$ is not $\Sigma_\text{cont}$--informative for stabilization. Therefore, $u_{[0,T-1]}$ is not $\Sigma_\text{cont}$--universal for stabilization, and thus, (b) does not hold. 

(c)$\Rightarrow$(a): Suppose that $u_{[0,T-1]}$ is persistently exciting of order $n+1$. Let $(A,B)\in\Sigma_\text{stab}$ and $x_{[0,T]}$ be such that $\mathcal{D}=(u_{[0,T-1]},x_{[0,T]})\in\mathfrak{B}_T(A,B)$. It follows from Lemma~\ref{lem:PE_stab} that both conditions in \eqref{eq:inf_stab_cond} hold. In case $\rank X_-<n$, $\mathcal{D}$ is $\Sigma_\text{stab}$--informative for stabilization due to Proposition~\ref{prop:inf_stab_pk(c)}.  Now, suppose that $\rank X_-=n$. In this case, $\rank\begin{bmatrix}X_-^\top & U_-^\top
\end{bmatrix}=n+m$. Therefore, due to Theorem~\ref{th:inf_id_pk}, $\mathcal{D}$ is $\mathcal{M}$--informative for identification, i.e., $\Sigma_\mathcal{D}=\{(A,B)\}$. Since $(A,B)$ is stabilizable, $\mathcal{D}$ is $\Sigma_\text{stab}$--informative for stabilization. Therefore, $u_{[0,T-1]}$ is $\Sigma_\text{stab}$--universal for stabilization.  \hfill $\blacksquare$
\end{pf}

Recall from Corollary~\ref{cor:cont_is_largest} that $\Sigma_\text{cont}$ is the largest set for which there exists a universal input for \emph{identification}. Now, due to Theorem~\ref{th:2}, there exist $\Sigma_\text{stab}$--universal inputs for stabilization. Hence, unlike identification, universal experiment design for stabilization can be feasible for sets of prior knowledge $\Sigma_\text{pk}$ that are larger than $\Sigma_\text{cont}$. On the other hand, it is easy to see that if $\Sigma_\text{pk}$ includes an element that is not stabilizable, then there are no $\Sigma_\text{pk}$--universal inputs for stabilization. Therefore, $\Sigma_\text{stab}$ is the largest set for which there exists a universal input for \emph{stabilization}. This is an interesting observation in parallel to the one made in Corollary~\ref{cor:cont_is_largest}.

\begin{corollary}
\label{cor:stab_is_largest}
There exists a $\Sigma_\text{pk}$--universal input for stabilization \emph{if and only if} $\Sigma_\text{pk}\subseteq\Sigma_\text{stab}$. 
\end{corollary}

\section{Online Experiment Design}
\label{sec:online}

In this section, we consider \emph{online} experiment design for identification and stabilization. In this approach, one selects the input at each step based on the data collected from the previous steps. This is in contrast to the offline nature of universal experiment design discussed in the previous sections. Therefore, one expects that an online approach leads to a shorter experiment compared to the design using a universal input. Here, we consider the online experiment design of Algorithm~\ref{alg:online}. Lemma~\ref{lem:online} will further elaborate on the significance of this online approach. 

\begin{algorithm}
\caption{Online experiment design}
\begin{algorithmic}[1]
\Require $x(0)\in\mathbb{R}^n$
\State set $t\gets 0$ and select $u(0)\neq 0$
\While{true}
\State apply $u(t)$ to $(A_\text{true},B_\text{true})$ and measure $x(t+1)$
\State set $t\gets t+1$
\If{$x(t)\notin \im X_{[0,t-1]}$}
    \State select $u(t)$ arbitrarily
\Else
    \If{$\im \begin{bmatrix}
    X_{[0,t-1]} \\
    U_{[0,t-1]}
    \end{bmatrix}= \im X_{[0,t-1]}\times \mathbb{R}^m$}
    \State \textbf{break}
    \Else
    \State select $\xi\in\mathbb{R}^n$ and $\eta\in\mathbb{R}^m\backslash\{0\}$ such that $$\xi^\top X_{[0,t-1]}+\eta^\top U_{[0,t-1]}=0$$
    \State select $u(t)$ such that $\xi^\top x(t)+\eta^\top u(t)\neq 0$
    \EndIf
\EndIf
\EndWhile
\Ensure $T=t$ and $\mathcal{D}=(u_{[0,T-1]},x_{[0,T]})$
\end{algorithmic}
\label{alg:online}
\end{algorithm}

\begin{lemma}
\label{lem:online}
Let $x(0)\in\mathbb{R}^n$ and $\mathcal{D}=(u_{[0,T-1]},x_{[0,T]})$ be generated by Algorithm~\ref{alg:online}. Then:
\begin{enumerate}[label=(\alph*),ref=\ref{lem:online}(\alph*)]
    \item\label{lem:online(a_true)} $T=\dim\mathcal{R}(A_\text{true},\begin{bmatrix}
    B_\text{true} & x(0)
    \end{bmatrix})+m$.
    \item\label{lem:online(a)} If $(A_\text{true},B_\text{true})\in\Sigma_\text{cont}$, then $\mathcal{D}$ is $\Sigma_\text{cont}$--informative for identification. 
    \item\label{lem:online(b)} If $(A_\text{true},B_\text{true})\in\Sigma_\text{stab}$, then $\mathcal{D}$ is $\Sigma_\text{stab}$--informative for stabilization. 
\end{enumerate}
\end{lemma}
\begin{pf}
(a) To simplify the notation, denote $T_\text{true}\coloneqq \dim\mathcal{R}(A_\text{true},\begin{bmatrix}
B_\text{true} & x(0)
\end{bmatrix})+m$. 
We show that the data $\mathcal{D}=(u_{[0,T-1]},x_{[0,T-1]})$ satisfies
\begin{equation}
\label{eq:lem:online-1}
\rank \begin{bmatrix}
X_{[0,t-1]} \\
U_{[0,t-1]}
\end{bmatrix}<\rank \begin{bmatrix}
X_{[0,t]} \\
U_{[0,t]}
\end{bmatrix}
\end{equation}
for all $t\in[1,T_\text{true}-1]$. For this, let $t\in[1,T_\text{true}-1]$. If $x(t)\notin \im X_{[0,t-1]}$, then \eqref{eq:lem:online-1} obviously holds. Now, suppose that $x(t)\in \im X_{[0,t-1]}$. We claim that 
\begin{equation}
\label{eq:condition_subspaces}
\im \begin{bmatrix}
    X_{[0,t-1]} \\
    U_{[0,t-1]}
    \end{bmatrix}=\im X_{[0,t-1]}\times \mathbb{R}^m
\end{equation}
does not hold. To show this, assume on the contrary that \eqref{eq:condition_subspaces} holds. Multiply \eqref{eq:condition_subspaces} from left by $\begin{bmatrix}
A_\text{true} & B_\text{true}
\end{bmatrix}$ to have $\im X_{[1,t]}=A_\text{true}\im X_{[0,t-1]} + \im B_\text{true}$. Since $x(t)\in \im X_{[0,t-1]}$, we have $\im X_{[1,t]}\subseteq\im X_{[0,t-1]}$. This implies that $\im X_{[0,t-1]}$ is $A_\text{true}$--invariant and contains $\im \begin{bmatrix}
B_\text{true} & x(0)
\end{bmatrix}$. Thus, we have $\mathcal{R}(A_\text{true},\begin{bmatrix}
    B_\text{true} & x(0)
    \end{bmatrix})\subseteq \im X_{[0,t-1]}$. This implies that $t\geq T_\text{true}$, which contradicts $t\leq T_\text{true}-1$. Therefore, \eqref{eq:condition_subspaces} does not hold. Hence, there exist $\xi\in\mathbb{R}^n$ and nonzero $\eta\in\mathbb{R}^m$ satisfying $\xi^\top X_{[0,t-1]}+\eta^\top U_{[0,t-1]}=0$. It is then also clear that there exists $u(t)$ such that $\xi^\top x(t)+\eta^\top u(t)\neq 0$. For this $u(t)$, we have that the dimension of $\ker\begin{bmatrix}
X_{[0,t]}^\top &
U_{[0,t]}^\top
\end{bmatrix}$ is strictly less than that of $\ker\begin{bmatrix}
X_{[0,t-1]}^\top &
U_{[0,t-1]}^\top
\end{bmatrix}$. Therefore, it follows from the rank-nullity theorem that \eqref{eq:lem:online-1} holds. Now, to show that $T=T_\text{true}$, we observe that $\im X_{[0,T_\text{true}-1]}\subseteq\im\mathcal{R}(A_\text{true},\begin{bmatrix}
B_\text{true} & x(0)
\end{bmatrix})$, and thus, we have
\begin{equation}
\label{eq:lem:online-2}
\im \begin{bmatrix}
X_{[0,T_\text{true}-1]} \\
U_{[0,T_\text{true}-1]}
\end{bmatrix}\subseteq \im\mathcal{R}(A_\text{true},\begin{bmatrix}
B_\text{true} & x(0)
\end{bmatrix})\times \mathbb{R}^m.
\end{equation}
It is also clear from \eqref{eq:lem:online-1} that  $\rank\begin{bmatrix}
X_{[0,T_\text{true}-1]} \\
U_{[0,T_\text{true}-1]}
\end{bmatrix}=T_\text{true}$. This, together with \eqref{eq:lem:online-2}, implies that
\begin{equation}
\label{eq:lem:online-3}
\im \begin{bmatrix}
X_{[0,T_\text{true}-1]} \\
U_{[0,T_\text{true}-1]}
\end{bmatrix}= \im\mathcal{R}(A_\text{true},\begin{bmatrix}
B_\text{true} & x(0)
\end{bmatrix})\times \mathbb{R}^m,
\end{equation}
and hence, 
\begin{equation}
\label{eq:lem:online-4}
\im X_{[0,T_\text{true}-1]}=\im\mathcal{R}(A_\text{true},\begin{bmatrix}
B_\text{true} & x(0)
\end{bmatrix}).
\end{equation}
In addition, since
\begin{equation}
\label{eq:lem:online-5}
x(T_\text{true})\in\im\mathcal{R}(A_\text{true},\begin{bmatrix}
B_\text{true} & x(0)
\end{bmatrix}),
\end{equation}
it is evident that at $t=T_\text{true}$ Algorithm~\ref{alg:online} terminates. Therefore, we have $T=T_\text{true}$.

(b) If $(A_\text{true},B_\text{true})$ is controllable, then we have that $\mathcal{R}(A_\text{true},\begin{bmatrix}
B_\text{true} & x(0)
\end{bmatrix})=\mathbb{R}^n$ for all $x(0)\in\mathbb{R}^n$. Therefore, it follows from statement (a) and the rank inequality \eqref{eq:lem:online-1} that $T=T_\text{true}=n+m$ and $\mathcal{D}$ satisfies $\rank \begin{bmatrix}
X_-^\top & U_-^\top
\end{bmatrix}=n+m$. Therefore, based on Theorem~\ref{th:inf_id_pk(b)} and Corollary~\ref{cor:inf_id(b)} we conclude that $\mathcal{D}$ is $\Sigma_\text{cont}$--informative for identification. 

(c) Due to \eqref{eq:lem:online-3} and  \eqref{eq:lem:online-4}, one can see that $\mathcal{D}$ satisfies $\im\begin{bmatrix}
    X_-^\top & U_-^\top
    \end{bmatrix}^\top=\im X_- \times \mathbb{R}^m$. In case $\rank X_-=n$, we have  $\rank\begin{bmatrix}
    X_-^\top & U_-^\top
    \end{bmatrix}=n+m$, and therefore, $\mathcal{D}$ is $\mathcal{M}$--informative for identification. This, in turn, implies that $\mathcal{D}$ is $\Sigma_\text{stab}$--informative for stabilization. In case $\rank X_-<n$, we observe that \eqref{eq:lem:online-4} and \eqref{eq:lem:online-5} imply $\im X_+\subseteq\im X_-$. Therefore, due to Proposition~\ref{prop:inf_stab_pk(c)}, $\mathcal{D}$ is $\Sigma_\text{stab}$--informative for stabilization. \hfill $\blacksquare$
\end{pf}

In case the true system is controllable, Algorithm~\ref{alg:online} boils down to the one presented in \citep{van2021beyond} and the stopping criterion reduces to checking whether $t=n+m$. In that case, it follows
from Lemma~\ref{lem:online(a_true)} that $T = n + m$, thus, $\mathcal{D}$ is the shortest dataset for identification. On the other hand, in case the true system is not controllable, collecting informative data for identification might not be possible. Interestingly, Lemma~\ref{lem:online(b)} shows that, in case the true system is \emph{stabilizable}, Algorithm~\ref{alg:online} generates $\Sigma_\text{stab}$--informative data for \emph{stabilization}. A similar observation was also made by \cite{gramlich2023fast}. 

Now, we investigate whether Algorithm~\ref{alg:online} generates the shortest experiment for stabilization. In order to facilitate the expressions, we introduce some notation. We denote the set of all $T$--length trajectories of the true system that are $\Sigma_\text{stab}$--informative for stabilization by 
\begin{equation}
\begin{split}
\mathfrak{D}_T\coloneqq \{\mathcal{D}\in\mathfrak{B}_T(A_\text{true},B_\text{true})\mid\mathcal{D}\text{ is } \Sigma_\text{stab}\text{--informative} \\ \text{for stabilization}\}.
\end{split}
\end{equation}
Based on Proposition \ref{prop:inf_stab_pk}, we know that if $T<\min\{n,m\}$, then $\mathfrak{D}_T=\varnothing$. Here, we are interested in the smallest value of $T\in\mathbb{N}$ such that $\mathfrak{D}_T$ is nonempty, i.e., the length of the shortest experiment. To formalize this, we define the following quantity as a function of the initial state $x_0\in\mathbb{R}^n$:
\begin{equation}
\begin{split}
T_*(x_0)\coloneqq\min\{T\in\mathbb{N}\mid\exists (u_{[0,T-1]},x_{[0,T]})\in\mathfrak{D}_T\\
\text{ with }x(0)=x_0\}.
\end{split}
\end{equation}
The value of $T_*(x_0)$ is equal to the number of samples of the shortest dataset that is $\Sigma_\text{stab}$--informative for stabilization starting from $x(0)=x_0$. Assuming that $(A_\text{true},B_\text{true})$ is stabilizable, it follows from Lemma~\ref{lem:online} that $T_*(x_0)$ satisfies $\min\{n,m\}\leq T_*(x_0)\leq n+m$. 

But how can we generate $\Sigma_\text{stab}$--informative data for stabilization with the length of the data exactly equal to $T_*(x_0)$? To the best of the authors' knowledge, there are no answers to this question in the literature. Here, we provide a partial answer to this question by showing that \emph{for some} values of $x_0$, the data generated by Algorithm~\ref{alg:online} will have a length equal to $T_*(x_0)$. Consider those $x_0$ satisfying
\begin{equation}
\label{eq:adverserial_x0}
\mathcal{R}(A_\text{true},\begin{bmatrix}
    B_\text{true} & x_0
    \end{bmatrix})\neq \mathbb{R}^n.
\end{equation}
If $x(0)=x_0$, then any dataset collected from the true system is not $\mathcal{M}$--informative for identification, no matter what the input signal is. Hence, these initial states are adversarial to system identification. Now, the following theorem shows that for such initial conditions, the data generated by Algorithm~\ref{alg:online} is surprisingly of the shortest possible length for stabilization. 

\begin{theorem}
\label{th:online_shortest}
Suppose that $(A_\text{true},B_\text{true})$ is stabilizable but not controllable. Let $\mathcal{D}$ be generated by Algorithm~\ref{alg:online} with $x(0)=x_0$, where $x_0$ satisfies \eqref{eq:adverserial_x0}. Then, the length of the data satisfies 
\begin{equation}
\label{eq:online_shortest}
T=T_*(x_0)= \dim\mathcal{R}(A_\text{true},\begin{bmatrix}
    B_\text{true} & x_0
    \end{bmatrix})+m.
\end{equation}
\end{theorem}
\begin{pf}
Since $\mathcal{D}$ is generated by Algorithm~\ref{alg:online}, it follows from Lemma~\ref{lem:online(a_true)} that its length satisfies $T=T_\text{true}\coloneqq\dim\mathcal{R}(A_\text{true},\begin{bmatrix}B_\text{true} & x_0\end{bmatrix})+m$. Moreover, in view of \linebreak Lemma~\ref{lem:online(a)}, $\mathcal{D}$ is $\Sigma_\text{stab}$--informative for stabilization. Thus, we have $T_*(x_0)\leq T=T_\text{true}$. What remains to be proven is that if $x_0$ satisfies \eqref{eq:adverserial_x0}, then $T_*(x_0)\geq T_\text{true}$. Let $\mathcal{D}^\prime=(u_{[0,k-1]},x_{[0,k]})\in\mathfrak{D}_k$ such that $x(0)=x_0$ satisfies \eqref{eq:adverserial_x0} and $k= T_*(x_0)$. One can verify that $\rank X_-<n$. Recall from Proposition~\ref{prop:inf_stab_pk(c)} that $\mathcal{D}^\prime$ satisfies $\im X_+\subseteq \im X_-$ and $\im\begin{bmatrix}X_-^\top & U_-^\top
\end{bmatrix}^\top=\im X_- \times \mathbb{R}^m$. Now, it follows from Lemma~\ref{lem:PE_stab} that $\im\begin{bmatrix}
X_-^\top & U_-^\top
\end{bmatrix}^\top=\mathcal{R}(A_\text{true},\begin{bmatrix}
B_\text{true} & x_0
\end{bmatrix}) \times \mathbb{R}^m$. This implies that $T_*(x_0)\geq T_\text{true}$. Therefore, $T_*(x_0)=T_\text{true}$, and thus, \eqref{eq:online_shortest} holds. \hfill $\blacksquare$
\end{pf}

\section{Conclusions and Discussion}

In Sections~\ref{sec:univ_id} and~\ref{sec:univ_stab}, we have presented necessary and sufficient conditions for an input to be universal for identification and stabilization while incorporating stabilizability or controllability as prior knowledge. It has been shown that persistency of excitation of order $n+1$ plays a central role in such characterizations. It has also been shown that prior knowledge on such system-theoretic properties is minimal in the sense that if the prior knowledge set is \emph{not} a subset of controllable (resp., stabilizable) systems, then experiment design for identification (resp., stabilization) is \emph{not} possible in general. 

In Section~\ref{sec:online}, we have extended our study to online experiment design using such prior knowledge. Unlike universal inputs that work simultaneously for \emph{all} systems within the prior knowledge set, an online approach is tailored towards a \emph{single} system. Each of these methods has its own pros and cons; while universal inputs are computationally inexpensive, online approaches are sample efficient. 

For both universal and online approaches, it turns out that experiment design for identification requires the prior knowledge set to be a subset of controllable systems. However, for stabilization, the prior knowledge needs to be a subset of stabilizable systems. This is an important observation, especially for systems that are not controllable but stabilizable; one can generate data suitable for stabilization when experiment design for identification is not possible.

Moreover, it has been shown that an online experiment design can generate the shortest dataset for stabilization, provided that the system is stabilizable but not controllable and the initial condition satisfies~\eqref{eq:adverserial_x0}. Interestingly, one can see from~\eqref{eq:online_shortest} that, e.g., the zero initial condition leads to the smallest number of samples. 

Future work may focus on other types of prior knowledge and the case where the data is corrupted by noise.

\bibliography{ifacconf}

@article{willems2005note,
  title={A note on persistency of excitation},
  author={Willems, Jan C and Rapisarda, Paolo and Markovsky, Ivan and De Moor, Bart LM},
  journal={Systems \& Control Letters},
  volume={54},
  number={4},
  pages={325--329},
  year={2005},
  publisher={Elsevier}
}

@article{van2021beyond,
  title={Beyond persistent excitation: Online experiment design for data-driven modeling and control},
  author={van Waarde, Henk J},
  journal={IEEE Control Systems Letters},
  volume={6},
  pages={319--324},
  year={2021},
  publisher={IEEE}
}

@article{camlibel2024beyond,
  title={Beyond the fundamental lemma: From finite time series to linear system},
  author={Camlibel, M. Kanat and Rapisarda, Paolo},
  journal={arXiv preprint arXiv:2405.18962},
  year={2024}
}

@article{van2020data,
  title={Data informativity: A new perspective on data-driven analysis and control},
  author={van Waarde, Henk J and Eising, Jaap and Trentelman, Harry L and Camlibel, M Kanat},
  journal={IEEE Transactions on Automatic Control},
  volume={65},
  number={11},
  pages={4753--4768},
  year={2020},
  publisher={IEEE}
}

@inproceedings{yu2021controllability,
  title={On controllability and persistency of excitation in data-driven control: Extensions of {W}illems’ fundamental lemma},
  author={Yu, Yue and Talebi, Shahriar and van Waarde, Henk J and Topcu, Ufuk and Mesbahi, Mehran and A{\c{c}}{\i}kmeșe, Beh{\c{c}}et},
  booktitle={IEEE Conference on Decision and Control},
  pages={6485--6490},
  year={2021}
}

@article{camlibel2025shortest,
  title={The shortest experiment for linear system identification},
  author={Camlibel, M K and van Waarde, H J and Rapisarda, P},
  journal={Systems \& Control Letters},
  volume={197},
  pages={106045},
  year={2025},
  publisher={ELSEVIER SCIENCE BV}
}

@article{shakouri2025new,
  title={A New Perspective on {W}illems’ Fundamental Lemma: Universality of Persistently Exciting Inputs},
  author={Shakouri, Amir and van Waarde, Henk J and Camlibel, M Kanat},
  journal={IEEE Control Systems Letters},
  volume={9},
  pages={583--588},
  year={2025},
  publisher={IEEE}
}

@article{shakouri2025pk,
  title={Data-Driven Stabilization Using Prior Knowledge on Stabilizability and Controllability},
  author={Shakouri, Amir and van Waarde, Henk J and Baltussen, Tren M. J. T. and  Heemels, W. P. M. H.},
  journal={arXiv preprint arXiv:2510.25452},
  year={2025}
}

@article{gramlich2023fast,
  title={Fast identification and stabilization of unknown linear systems},
  author={Gramlich, Dennis and Ebenbauer, Christian},
  journal={IFAC-PapersOnLine},
  volume={56},
  number={2},
  pages={6241--6246},
  year={2023},
  publisher={Elsevier}
}

@article{berberich2025overview,
  title={An overview of systems-theoretic guarantees in data-driven model predictive control},
  author={Berberich, Julian and Allg{\"o}wer, Frank},
  journal={Annual Review of Control, Robotics, and Autonomous Systems},
  volume={8},
  number={1},
  pages={77--100},
  year={2025},
  publisher={Annual Reviews}
}

@book{van2025data,
  title={Data-Based Linear Systems and Control Theory},
  author={van Waarde, Henk J and Camlibel, M Kanat and Trentelman, Harry L},
  year={2025},
  publisher={Kindle Direct Publishing}
}

@article{markovsky2021behavioral,
  title={Behavioral systems theory in data-driven analysis, signal processing, and control},
  author={Markovsky, Ivan and D{\"o}rfler, Florian},
  journal={Annual Reviews in Control},
  volume={52},
  pages={42--64},
  year={2021},
  publisher={Elsevier}
}

\end{document}